\documentclass[11pt,a4paper]{article}

\usepackage{epsf,epsfig,amsfonts,amsgen,amsmath,amstext,amsbsy,amsopn,amsthm
}
\usepackage{amsmath,times,mathptmx}
\usepackage{amsfonts,amsthm,amssymb}
\usepackage{amsfonts}
\usepackage{graphics}
\usepackage{latexsym,bm}
\usepackage{amsfonts,amsthm,amssymb,bbding}
\usepackage{indentfirst}
\usepackage{graphicx}
\usepackage{color}
\usepackage[colorlinks=true,anchorcolor=blue,filecolor=blue,linkcolor=blue,urlcolor=blue,citecolor=blue]{hyperref}
\usepackage{float}
\usepackage{tikz}
\newtheorem{thm}{Theorem}[section]

\newtheorem{lemma}{Lemma}[section]
\newtheorem{prop}[lemma] {Proposition} 
\newtheorem{false statement}{False statement}
\newtheorem{cor}{Corollary} [section]

\theoremstyle{definition}

\newtheorem{claim}{Claim}

\begin{document}

\title{Spectral radius and rainbow $k$-factors of graphs
}
\author{{Liwen Zhang, Zhiyuan Zhang}\thanks{Corresponding author. E-mail addresses: jsam0331@163.com (Z. Zhang), levenzhang512@163.com (L. Zhang).}\\
{\footnotesize Center for Combinatorics and LPMC, Nankai University, Tianjin 300071, China}}
\date{}

\maketitle{\flushleft\large\bf Abstract:}
 Let $\mathcal{G}=\{G_1,\ldots, G_{\frac{kn}{2}}\}$ be a set of graphs on the same vertex set $V=\{1,\dots,n\}$ where $k\cdot n$ is even. We say $\mathcal{G}$ admits a rainbow $k$-factor if there exists a $k$-regular graph $F$ on the vertex set $V$ such that all edges of $F$ are from different members of $\mathcal{G}$. In this paper, we show a sufficient spectral condition for the existence of a rainbow $k$-factor for $k\geq 2$, which is that if $\rho(G_i)\geq\rho(K_{k-1}\vee(K_1\cup K_{n-k}))$ for each $G_i\in \mathcal{G}$, then $\mathcal{G}$ admits a rainbow $k$-factor unless $G_1=G_2=\cdots=G_{\frac{kn}{2}}\cong K_{k-1}\vee(K_1\cup K_{n-k})$.

\vspace{0.1cm}
\begin{flushleft}
\end{flushleft}

\section{Introduction}
The spanning subgraphs of graphs possessing some given properties are called \textit{factors}. Factor theory is one of the fundamental areas of graph theory and closely relates to graph spectra. A considerable number of influential results on spectra and factors have been established, see, for example, the survey \cite{FLLO23}. This paper primarily focuses on \textit{$[a, b]$-factor}, which is defined as a spanning subgraph $H$ of $G$ such that $d_H(v)\in [a,b]$ for each $v\in V (G)$. Particularly, a $[k,k]$-factor is also called a \textit{$k$-factor}. 

In 2005, Brouwer and Haemers \cite{BH2005} provided a sufficient condition for the existence of a perfect matching (1-factor) in a regular graph in terms of the third largest eigenvalue. Subsequently, this condition was improved in \cite{Cioaba2005,CG2007,CGH2009} and extended to $k$-factor in \cite{Lu2010,Lu2012}. In terms of spectral radius, Fiedler and Nikiforov \cite{FN10} showed tight conditions on the spectral radius that guarantee the existence of a Hamiltonian path (a special $[1,2]$-factor) and a Hamiltonian cycle (a special $2$-factor), respectively. In 2021, O \cite{O2021} gave an upper bound on the spectral radius for graphs that do not contain perfect matchings. Denote by $\vee$ and $\cup$ the \emph{join} and \emph{union} operations on graphs, respectively. Let $H_{n,k}\cong K_{k-1}\vee(K_1\cup K_{n-k})$. Fan, Lin and Lu \cite{FLL22} showed that if $\rho(G)>\rho(H_{n,a})$ and $n\geq 3a+b+1$, then the graph $G$ contains an $[a,b]$-factor, which 
partially confirmed the conjecture on $[a,b]$-factor proposed by Cho, Hyun, O and Park \cite{CHOP2021}. Further, their result also gave the spectral condition for the existence of a $k$-factor.
\begin{thm}[\cite{FLL22}, Corollary 1.4]\label{kfactor}
    Let $k$ and $n\geq 4k-1$ be two positive integers such that $k\cdot n$ is even and let $G$ be a graph of order $n$. If $\rho(G)>\rho(H_{n,k})$, then $G$ contains a $k$-factor.
\end{thm} 

 Let $\mathcal{G}=\{G_1,G_2,\dots,G_m\}$ be a collection of (not necessarily distinct) graphs on the same vertex set $V$, and let $H$ be a graph with $m$ edges on the vertex set $V(H)\subseteq V$. If there exists a bijection $\varphi: E(H)\to \{1,\ldots,m\}$ such that $e\in E(G_{\varphi(e)})$ for each $e\in E(H)$, then we say that $H$ is a \textit{rainbow} of $\mathcal{G}$, or $\mathcal{G}$ admits a \textit{rainbow} $H$. Recently, the study of rainbow structures has attracted extensive attention, yielding many beautiful results from a variety of perspectives, including degree \cite{B21,CSWW24+,GHMPS23,JK20,LLL23}, size \cite{ADGM20,HFGLSTKZ24,KSSV04,M15} and spectra \cite{GLMM23,HLF24,ZV24+}. For more results, we refer the readers to the survey \cite{SWW24+}.
 
In this paper, we discuss the spectral condition on the set of graphs $\mathcal{G}=\{G_1,G_2,\dots,G_m\}$ that guarantee the existence of a rainbow factor in $\mathcal{G}$. Guo, Lu, Ma, and Ma \cite{GLMM23} identified the spectral conditions related to rainbow matchings. Later, He, Li and Feng \cite{HLF24} gave sufficient spectral conditions for rainbow Hamiltonian paths and rainbow linear forests. Recently, Zhang and van Dam \cite{ZV24+} proved that $\mathcal{G}=\{G_1, G_2,\cdots, G_{n}\}$ admits a rainbow Hamiltonian cycle if $\rho(G_i)>n-2$ for $i=1,2,\ldots,n$ unless $G_1=G_2=\cdots=G_{n}\cong H_{n,2}$. Inspired by these extensive works, we consider an extension of Theorem \ref{kfactor} to the rainbow
version, by presenting a spectral condition for a rainbow $k$-factor. 
\begin{thm}\label{rainbowkfactor}
    Let $k$ and $n\geq 4k-1$ be two positive integers such that $k\cdot n$ is even and let $\mathcal{G}=\{G_1,\ldots, G_{\frac{kn}{2}}\}$ be a set of graphs on the same vertex set $V=\{1,\dots,n\}$. If $\rho(G_i)\geq\rho(H_{n,k})$ for $i=1,2,\ldots,\frac{kn}{2}$, then $\mathcal{G}$ admits a rainbow $k$-factor unless $G_1=G_2=\cdots=G_{\frac{kn}{2}}\cong H_{n,k}$.
\end{thm}

Our paper is organized as follows: In Section 2, we introduce the so-called Kelmans operation and do some preliminary
work. In Section 3, we provide the proof of our main result. In addition, we refer the readers to \cite{BM76,BH12} for notation and terminology not defined here.

\section{Preliminaries}
 Let $G$ be a simple graph with vertex set $V(G)=\{1,\ldots,n\}$ and edge set $E(G)$.  The \textit{adjacency matrix} of $G$ is defined as $A(G)=(a_{ij})_{n\times n}$ where $a_{ij}=1$ if $ij\in E(G)$ and 0 otherwise. The \textit{eigenvalues} of $G$ are the eigenvalues of $A(G)$. The largest eigenvalue is also known as the
 \textit{spectral radius} of $G$, denoted by $\rho(G)$. Two graphs $G_1$ and $G_2$ are said to be \textit{isomorphic}, denoted $G_1\cong G_2$, if there exists a bijection $\phi:V (G_1)\rightarrow V(G_2)$ such that $uv\in E(G_1)$ if and only if $\phi (u)\phi(v)\in E(G_2)$. Furthermore, two graphs $G_1$ and $G_2$ are \textit{identical}, written $G_1=G_2$, if $V(G_1)=V(G_2)$ and $E(G_1)=E(G_2)$. 
 
Given a graph $G$ and $u,v\in V(G)$, we can obtain the graph ${KO}_{uv}(G)$ from $G$ via replacing the edge $vw$ by a new edge $uw$ for all $w\in N_G(v)\setminus (N_G(u)\cup \{u\})$. This operation is called the \textit{Kelmans operation} (from $v$ to $u$) due to Kelmans \cite{Kelmans1981}. In 2009, Csikv\'{a}ri \cite{C09} proved that the spectral radius of a graph does not decrease after the Kelmans operation.

\begin{lemma}[\cite{C09}, Theorem 2.1]\label{boundKO}
    Let $G$ be a graph on vertex set $V(G)$ and let $u,v\in V(G)$. Then $$\rho({KO}_{uv}(G))\geq\rho(G).$$
\end{lemma}

In 2019, Liu, Lai and Das found that the above inequality is strict when $G$ is connected and $G\not\cong {KO}_{uv}(G)$.

\begin{lemma}[\cite{LLD19}, Corollary 4.1]\label{newbound}
    Let $G$ be a connected graph on vertex set $V(G)$ and let $u,v\in V(G)$. Then $\rho({KO}_{uv}(G))>\rho(G)$ unless $G\cong {KO}_{uv}(G)$.
\end{lemma}

Observe that ${KO}_{uv}(G)\cong {KO}_{vu}(G)$. Denote by $KO(G)$ the graph obtained from $G$ by applying the Kelmans operation to each vertex pair $(u,v)$ such that $u<v$. It implies that $KO_{uv}(KO(G))\cong KO(G)$ for any vertex pair $(u,v)$. Then, we obtain the following proposition in $KO(G)$, which will be useful in the subsequent proofs.

\begin{prop}\label{prop}
Let $G$ be a graph on vertex set $V=\{1,\ldots,n\}$ and let $ x,y\in  V$ such that $x<y$. If $xy\in E(KO(G))$, then $ij\in E(KO(G))$ for each integer $1\leq i\leq x$ and $i< j\leq y$.
\end{prop}

 For a given graph set $\mathcal{G}=\{G_1,\dots,G_{m}\}$, define $KO_{uv}(\mathcal{G})=\{KO_{uv}(G_1),\dots,KO_{uv}(G_{m})\}$ and $KO(\mathcal{G})=\{KO(G_1),\dots,KO(G_{m})\}$.

\begin{lemma}\label{structural}
   Let $\mathcal{G}=\{G_1,\ldots,G_{\frac{kn}{2}}\}$ be a set of graphs on the same vertex set $V=\{1,\ldots,n\}$. For any vertices $u,v \in V$, if ${KO}(\mathcal{G})$ admits a rainbow $k$-factor, then so does $\mathcal{G}$. 
\end{lemma}

\begin{proof} It suffices to show that for any vertices $u,v \in V$, if ${KO}_{uv}(\mathcal{G})$ admits a rainbow $k$-factor, then so does $\mathcal{G}$.
Assume that $F$ is a rainbow $k$-factor of ${KO}_{uv}(\mathcal{G})$ where $E(F)=\{e_1,e_2,\ldots,e_{\frac{kn}{2}}\}$ and $e_i\in E({KO}_{uv}(G_i))$ for $i\in [1,\frac{kn}{2}]$. If $e_i\in E(G_i)$ for each integer $i\in [1,\frac{kn}{2}]$, then $F$ is also a rainbow $k$-factor of $\mathcal{G}$. Otherwise, there exists an edge $e_r=uw\in E(KO_{uv}(G_r))$ but $e_r\notin E(G_r)$. According to the definition of the Kelmans operation, we have $vw\in E(G_r)$. If $vw\in E(F)$, then without loss of generality, assume $e_s=vw\in E({KO}_{uv}(G_s))$. This implies both $vw\in E(G_s)$ and $uw\in E(G_s)$. In this case, let $e_r=vw$ and $e_s=uw$. If $vw\notin E(F)$, there is a vertex $w'$ satisfying $e_t=vw'\in E(F)$ and $uw'\notin E(F)$ since $d_{F}(u)=d_{F}(v)=k$. Note that $e_t=vw'\in KO_{uv}(G_t)$, which implies that $vw'\in G_t$ and $uw'\in G_t$. Replace $uw$ by $vw$ and $vw'$ by $uw'$, so that $e_r=vw$ and $e_t=uw'$. Then we obtain $e_i\in E(G_i)$ for any $i\in [1,\frac{kn}{2}]$. Hence, $\mathcal{G}$ admits a rainbow $k$-factor.
\end{proof}

Moreover, the proofs also involve several established results.

\begin{lemma}[\cite{BH12}]\label{lemsubgraph}
If $H$ is a subgraph of a connected graph $G$, then $\rho(H)\le \rho(G)$ with equality if and only if $H\cong G$.
\end{lemma}

\begin{thm}[\cite{GLMM23}, Theorem 3]\label{rainbowmatching}
 Let $\mathcal{G}=\{G_1,\ldots, G_{\frac{n}{2}}\}$ be a set of graphs on the same vertex set $V=\{1,\dots,n\}$ where $n$ is even. If $\rho(G_i)\geq n-2$ for $i=1,2,\ldots, \frac{n}{2}$, then $\mathcal{G}$ admits a rainbow perfect matching unless $G_1=G_2=\cdots=G_{\frac{n}{2}}\cong H_{n,1}$.
 \end{thm}

 \begin{thm}[\cite{ZV24+}, Theorem 4.3]\label{rainbowhamiltoncycle}
 Let $\mathcal{G}=\{G_1,\ldots, G_{n}\}$ be a set of graphs on the same vertex set $V=\{1,\dots,n\}$ where $n\geq 4$. If $\rho(G_i)> n-2$ for $i=1,2,\ldots,n$, then $\mathcal{G}$ admits a rainbow Hamiltonian cycle unless $G_1=G_2=\cdots=G_{n}\cong H_{n,2}$.
 \end{thm}
 
\section{Proof of Theorem \ref{rainbowkfactor}}


\begin{lemma}\label{xyKO}
    Let $G$ be a graph on vertex set $V=\{1,\ldots,n\}$ and let $u,v\in V$ such that $u<v$. If $\rho(G)=\rho(H_{n,k})$ and ${KO}_{uv}(G)\cong H_{n,k}$, then $G\cong H_{n,k}$.
\end{lemma}

\begin{proof}
   The condition ${KO}_{uv}(G)\cong H_{n,k}$ implies that there is a vertex subset $A\subset V$ of size $n-1$ such that ${KO}_{uv}(G)[A]\cong K_{n-1}$. Since $n-1\leq |\{u,v\}\cup A|\leq n$, we have $1\leq |\{u,v\}\cap A|\leq 2$. Observe that $G\cong KO_{uv}(G)\cong H_{n,k}$ when $|\{u,v\}\cap A|=2$. Thus, assume $|\{u,v\}\cap A|=1$ in the remaining proof. If $G$ is disconnected, then $\rho(G)\leq \rho(K_{n-1})=n-2\leq \rho(H_{n,k})$, and the equality holds if and only if $k=1$. It follows that $G\cong  H_{n,1}$. If $G$ is connected, then by Lemma \ref{newbound}, we have $\rho(KO_{uv}(G))>\rho(G)$ unless $G\cong KO_{uv}(G)\cong H_{n,k}$. Combining this with $\rho(G)=\rho(H_{n,k})=\rho(KO_{uv}(G))$ yields $G\cong H_{n,k}$. Hence, the proof is complete.
\end{proof}

  From Lemma \ref{boundKO} and Lemma \ref{xyKO}, we can deduce the following corollary.

\begin{cor}\label{cor}
    Let $G$ be a graph of order $n$. If $\rho(G)=\rho(H_{n,k})$ and $KO(G)\cong H_{n,k}$, then $G\cong H_{n,k}$.
\end{cor}

Recall that the join and union operations are denoted by $\vee$ and $\cup$, respectively. Moreover, the union of disjoint sets is denoted by $\Dot{~\cup~}$. For simplicity, let $ [n]=\{1,\dots,n\}$.

\begin{lemma}\label{factlemma}
    Let $k\geq 2$ and $n\geq 4k-1$ be two integers such that $k\cdot n$ is even and let $\mathcal{G}=\{G_1, \ldots, G_{\frac{kn}{2}}\}$ be a set of graphs on the same vertex set $V=\{1,\ldots,n\}$. If $G_i\cong H_{n,k}$ for each $i\in [\frac{kn}{2}]$, and there exist $i_1,i_2\in [\frac{kn}{2}]$ such that $G_{i_1}\ne G_{i_2}$, then $\mathcal{G}$ admits a rainbow $k$-factor.
\end{lemma}

\begin{proof}
    For the simplicity, let $V(G_i)=A_i\Dot{~\cup~}B_i\Dot{~\cup~}C_i$ such that $G_i[A_i]\cong K_1$, $G_i[B_i]\cong K_{k-1}$, $G_i[C_i]\cong K_{n-k}$, and $G_i=G_i[B_i]\vee(G_i[A_i]\cup G_i[C_i])$ for $1\leq i\leq \frac{kn}{2}$. 
 \begin{claim}\label{claimequi}
        If $A_1=A_2=\cdots=A_{\frac{kn}{2}}$, then $\mathcal{G}$ admits a rainbow $k$-factor.
    \end{claim}
    \begin{proof}
Without loss of generality, let $A_1=A_2=\cdots=A_{\frac{kn}{2}}=\{u\}$ and $G_1\neq G_2$. Then $|\cup_{i=1}^k B_i|\geq k$. Thus, there exists a subset $B'=\{v_1,\ldots,v_k\}\subset V$ such that $v_i\in B_i$ for $1\leq i \leq k$. Let $A'=\{u\}$ and $C'=V\setminus(A'\cup B')$. Now we construct a graph $H$ on the vertex set $V$ that satisfies $H[B']\cong K_{k}$, $H[C']\cong K_{n-k-1}$ and $H=H[B']\vee(H[A']\cup H[C'])$. Note that $\rho(H)>\rho(H_{n,k})$. By Theorem \ref{kfactor}, the graph $H$ contains a $k$-factor $F$ with $E(F)=\{e_1,e_2,\ldots,e_{\frac{kn}{2}}\}$. As $d_{H}(u)=k=d_{F}(u)$, we obtain $\{uv_1,uv_2,\ldots,uv_k\}\subset E(F)$. So, assume $e_i=uv_i$ when $1\leq i\leq k$. The definition of $B'$ implies that $e_i\in E(G_i)$ for each $1\leq i\leq k$. Furthermore, for any $i\in \{k+1,\ldots,\frac{kn}{2}\}$, we have $e_i\in E(G_i)$. Therefore, $\mathcal{G}$ admits a rainbow $k$-factor.\end{proof}

  Next, partition $\mathcal{G}$ as $\mathcal{G}=\mathcal{G}_1\Dot{~\cup~}\mathcal{G}_2\Dot{~\cup~}\cdots\Dot{~\cup~}\mathcal{G}_q$ satisfies 
  \begin{itemize}
      \item [\rm 1)] for $i\in [1,q]$, $n_i=|\mathcal{G}_i|$ and $n_1\geq n_2\geq \dots\geq n_q$;
       \item [\rm 2)]  for any two graphs $G_i\in \mathcal{G}_{\ell}$ and $G_j\in \mathcal{G}_{m}$, if $\ell=m$, then $A_i=A_j$; otherwise $A_i\neq A_j$, where $i,j\in [\frac{kn}{2}]$ and $\ell,m\in [q]$.
  \end{itemize}
Note that $q\geq 2$ by Claim \ref{claimequi}. For any $G_i\in \mathcal{G}_\ell$, assume that $A_{i}=\{u_\ell \}$ where $1\leq \ell \leq q$. Then we divide the proof into two cases depending on whether $n$ is even or odd.   

  \noindent{{\underline{{\bf Case 1. }$n$ is even.}}} 

  Let $\hat{k}_i=\lfloor\frac{n_i}{n/2}\rfloor$ for $i\in [q]$. Note that $\hat{k}_1\geq \hat{k}_2 \geq \dots\geq \hat{k}_q\geq 0$. When $\hat{k}_1>0$, define $p\leq q$ as the largest integer such that $\hat{k}_p>0$ and let $\hat{k}_1+\hat{k}_2+\dots+\hat{k}_p=\hat{k}\leq k$; when $\hat{k}_1=0$, let $p=0$ and $\hat{k}=0$. Set $\mathcal{G}'=\mathcal{G}'_1\Dot{~\cup~}\mathcal{G}'_2\Dot{~\cup~}\cdots\Dot{~\cup~}\mathcal{G}'_p$ where $\mathcal{G}'_i\subseteq \mathcal{G}_i$ and $|\mathcal{G}'_i|=\frac{\hat{k}_i n}{2}$ for $i\in [p]$. Then $|\mathcal{G}'|=\frac{\hat{k} n}{2}$ and we have the following claim.

\begin{claim}\label{claimeven}
    The graph set $\mathcal{G}'$ admits a rainbow $\hat{k}$-factor.
    \end{claim}

    \begin{proof}
    The result is trivial for $p=0$. According to Claim \ref{claimequi}, it suffices to show that the result holds for $p\geq 2$. Assume that $\mathcal{G}'_\ell=\{G_{\ell,1},G_{\ell,2},\ldots,G_{\ell, \frac{\hat{k}_\ell n}{2}}\}$ for $\ell \in [1,p]$. Similarly as in the proof of Claim \ref{claimequi}, we deduce that each $\mathcal{G}'_\ell$ admits a rainbow $\hat{k}_\ell$-factor $F_\ell$ where $E(F_\ell)=\{e_{\ell,1},~e_{\ell,2},~\ldots,~e_{\ell, \frac{\hat{k}_\ell n}{2}}\}$ such that $e_{\ell,i}\in E(G_{\ell,i})$ for $1\leq i \leq \frac{\hat{k}_\ell n}{2}$. If $E(F_r) \cap E(F_s)=\emptyset$ for any $1\leq r<s\leq p$, then $F_1\cup F_2 \cup\cdots \cup F_p$ is a rainbow $\hat{k}$-factor of $\mathcal{G}'$. Otherwise, suppose that $e_{r,i_1}=e_{s,i_2}=vv'$. Note that $\hat{k}\geq 2$ by $p\geq 2$. Then we could find an edge $ww'=e_{t,i_3}\in E(F_t)$ satisfying $w,w'\in V\setminus\{v,v'\}$ and $vw,vw',v'w,v'w'\notin E(F)$, since
    \begin{equation*}
        \frac{\hat{k}n}{2}-2-2\cdot(\hat{k}-2)-2(\hat{k}-2)(\hat{k}-1)=\frac{\hat{k}n}{2}-2\hat{k}(\hat{k}-2)-2\geq 1.
    \end{equation*}
    Without loss of generality, let $r \neq t$. Notice that at least one of the vertices $v$ and $v'$ is not equal to $u_r$, and at least one is not equal to $u_t$. Assume $v\neq u_r$ and $v'\neq u_t$, similarly $w\neq u_t$ and $w'\neq u_r$. Then replace $vv'\in E(G_{r,i_1})$ by $vw'\in E(G_{r,i_1})$ and $ww'$ by $v'w$. It implies that $e_{r,i_1}=vw'$, $e_{s,i_2}=vv'$ and $e_{t,i_3}=v'w$ after the replacement. Repeating the above steps until $E(F_r) \cap E(F_s)=\emptyset$ for $1\leq r<s\leq p$, then we obtain a rainbow $\hat{k}$-factor of $\mathcal{G}'$.
    \end{proof}

From Claim \ref{claimeven}, let $F'$ be a rainbow $\hat{k}$-factor of $\mathcal{G}'$ where $E(F')=\cup_{i=1}^p\{e_{i,1},\ldots,e_{i,\frac{\hat{k}_i n}{2}}\}$ such that $e_{i,j}\in E(G_{i,j})$ and $G_{i,j}\in \mathcal{G}'_i$. If $\hat{k}= k$, we are done. For $0\leq \hat{k}\leq k-1$, assume that $\mathcal{G}''=\mathcal{G}\setminus\mathcal{G}'=\{G_{\frac{\hat{k}n}{2}+1},\dots,G_{\frac{kn}{2}}\}$. Observe that for any $\frac{n}{2}$ graphs in the set $\mathcal{G}''$, there exist at least two graphs, say $G_{a}$ and $G_{b}$, such that $A_{a}\ne A_{b}$. Let $G'_{j}=G_j[A_{j}]\cup (G_j[B_j]\vee G_j[C_j])\cong H_{n,1}$ for $j=\frac{\hat{k}n}{2}+1,\dots,\frac{kn}{2}$. Using Theorem \ref{rainbowmatching}, we obtain that $\{G'_{\frac{\hat{k}n}{2}+1}.\dots,G'_{\frac{kn}{2}}\}$ admits $k-\hat{k}$ rainbow perfect matchings, denoted by $M_{\hat{k}+1},\ldots,M_k$. Then $M_{\hat{k}+1},\ldots,M_k$ are $k-\hat{k}$ rainbow perfect matchings of $\mathcal{G}''$ as well since $G'_{i}$ is a spanning subgraph of $G_{i}$. Without loss of generality, assume that $E(M_i)=\{e_{i,1},\dots,e_{i,\frac{n}{2}}\}$ such that $e_{i,j}\in E(G'_{\frac{(i-1)n}{2}+j})\subset E(G_{\frac{(i-1)n}{2}+j})$ where $i=\hat{k}+1,\dots,k$ and $j=1,\ldots,\frac{n}{2}$.

If $E(M_i)\cap E(M_j)=\emptyset$ and $E(F')\cap E(M_j)=\emptyset$ for any integers $\hat{k}+1\leq i,j\leq k$, we are done as $F'\cup M_{\hat{k}+1}\cup\cdots\cup M_{k}$ is a rainbow $k$-factor of $\mathcal{G}$. Otherwise, suppose that $e_{r,i_1}=e_{s,i_2}=vv'$ where $\hat{k} <r<s\leq k$ or $1\leq r\leq p\leq \hat{k} <s\leq k$. So we have $e_{s,i_2}\in E(G_{\frac{(s-1)n}{2}+i_2})$. Assume that $G_{\frac{(s-1)n}{2}+i_2}\in \mathcal{G}_{\ell_0}$ where $\ell_0 \in [ q]$. Notice that $v,v'\neq u_{\ell_0}$ since $vv'\in E(G'_{\frac{(s-1)n}{2}+i_2})$. Then we could find an edge $ww'=e_{t,i_3}\in E(F'\cup M_{\hat{k}+1}\cup\dots\cup M_{k})$ satisfying that $w, w'\in V\setminus\{v,v',u_{\ell_0} \}$ and $vw,vw',v'w,v'w'\notin E(F'\cup M_{\hat{k}+1}\cup\dots\cup M_{k})$, since
    \begin{equation*}
        \frac{kn}{2}-2-2\cdot(k-2)-k-2(k-2)(k-1)=\frac{kn}{2}-k(2k-3)-2\geq 1.
    \end{equation*}
    If $e_{t,i_3}\in E(F')$, then $e_{t,i_3}\in E(G_{t,i_3})$ and $G_{t,i_3}\in \mathcal{G}_{t}$; otherwise $e_{t,i_3}\in E(M_{\hat{k}+1}\cup\dots\cup M_{k})$ and $e_{t,i_3}\in E(G_{\frac{(t-1)n}{2}+i_3})$, without loss of generality, we could still assume that $G_{\frac{(t-1)n}{2}+i_3}\in \mathcal{G}_{t}$ (it is possible that $t=\ell_0$). Thus, all vertices in $\{v,v',w,w'\}$ are distinct from $u_{\ell_0}$ and at least one of the vertices $w$ and $w'$ ($v$ and $v'$) is distinct from $u_t$, say $w\neq u_t$ ($v\neq u_t$). Then replace $vv'\in E(G_{\frac{(s-1)n}{2}+i_2}) $ and $ww'$ by $v'w'\in E(G_{\frac{(s-1)n}{2}+i_2})$ and $vw$, respectively. It implies that $e_{r,i_1}=vv'$, $e_{s,i_2}=v'w'$ and $e_{t,i_3}=vw$ after the replacement. Repeating the above steps, we finally obtain a rainbow $k$-factor of $\mathcal{G}$. Hence, we are done with Case 1. 

    \noindent{\underline{{\textbf{Case 2. } $n$ is odd. }}} 
 
    In this case, $k$ is even. Let $\tilde{k}_i=\lfloor\frac{n_i}{n}\rfloor$ for $i \in [q]$. Observe that $\tilde{k}_1\geq \tilde{k}_2 \geq \dots\geq \tilde{k}_q\geq 0$. When $\tilde{k}_1>0$, define $p\leq q$ as the largest integer such that $\tilde{k}_p>0$ and let $\tilde{k}_1+\tilde{k}_2+\dots+\tilde{k}_p=\tilde{k}\leq \frac{k}{2}$; when $\tilde{k}_1=0$, let $p=0$ and $\tilde{k}=0$. Set $\mathcal{G}^{*}=\mathcal{G}^{*}_1\Dot{~\cup~}\mathcal{G}^{*}_2\Dot{~\cup~}\cdots\Dot{~\cup~}\mathcal{G}^{*}_p=\cup_{i=1}^p\{G_{i,1},\ldots,G_{i,\tilde{k}_i n}\}$ where $\mathcal{G}^{*}_i\subseteq \mathcal{G}_i$ and $|\mathcal{G}^{*}_i|=\tilde{k}_i n$ for $i \in [p]$. Similarly to the proof of Claim \ref{claimeven}, we have the following claim.

\begin{claim}\label{claimodd}
       The graph set $\mathcal{G}^{*}$ admits a rainbow $2\tilde{k}$-factor.
    \end{claim}

    According to Claim \ref{claimodd}, it suffices to show that the result holds for $0\leq \tilde{k}\leq \frac{k}{2}-1$. Let $F^*$ be a rainbow $2\tilde{k}$-factor of $\mathcal{G}^*$ where $E(F^*)=\cup_{i=1}^p\{e_{i,1},\ldots,e_{i,\tilde{k}_i n}\}$ such that $e_{i,j}\in E(G_{i,j})$ and $G_{i,j}\in \mathcal{G}^*_i$. Observe that for any ${n}$ graphs in the graph set $\mathcal{G}^{**}=\mathcal{G}\setminus\mathcal{G}^*=\{G_{{\tilde{k}n}+1},\dots,G_{\frac{kn}{2}}\}$, there exist at least two graphs, say $G_{a}$ and $G_{b}$, such that $A_{a}\neq A_{b}$. Using Theorem \ref{rainbowhamiltoncycle}, we deduce that $\mathcal{G}^{**}$ admits $\frac{k}{2}-\tilde{k}$ rainbow Hamiltonian cycles, denoted by $M_{\tilde{k}+1},\ldots,M_{\frac{k}{2}}$, where $E(M_i)=\{e_{i,1},\dots,e_{i,n}\}$ satisfies $e_{i,j}\in G_{(i-1)n+j}$ for $i=\tilde{k}+1,\ldots,\frac{k}{2}$ and $j=1,2,\dots,{n}$. 

    If $E(M_i)\cap E(M_j)=\emptyset$ and $E(F^*)\cap E(M_j)=\emptyset$ for any integers $\tilde{k}+1\leq i,j\leq \frac{k}{2}$, we are done as $\mathcal{G}$ admits a rainbow $k$-factor $F^*\cup M_{\tilde{k}+1}\cup\cdots\cup M_{\frac{k}{2}}$. Otherwise, suppose that $e_{r,i_1}=e_{s,i_2}=vv'$ where $\tilde{k}<r<s\leq \frac{k}{2}$ or $1\leq r\leq p\leq \tilde{k}<s\leq \frac{k}{2}$. Note that $e_{r,i_1}\in E(G_{(r-1)n+i_1})$ (or $e_{r,i_1}\in E(G_{r,i_1})$) and $e_{s,i_2}\in E(G_{(s-1)n+i_2})$. Without loss of generality, assume that $G_{(r-1)n+i_1}\in \mathcal{G}_{\ell_1}$ (or $G_{r,i_1}\in \mathcal{G}_{\ell_1}$) and $G_{(s-1)n+i_2}\in \mathcal{G}_{\ell_2}$ where $1\leq \ell_1,\ell_2\leq q$. Similarly, there exists an edge $ww'=e_{t,i_3}\in E(F^*\cup M_{\tilde{k}+1}\cup\dots\cup M_{\frac{k}{2}})$ such that $w, w'\in V\setminus\{v,v',u_{\ell_2} \}$ and $vw,vw',v'w,v'w'\notin E(F^*\cup M_{\tilde{k}+1}\cup\dots\cup M_{\frac{k}{2}})$. If $e_{t,i_3}\in E(F^*)$, we have $e_{t,i_3}\in E(G_{t,i_3})$ and $G_{t,i_3}\in \mathcal{G}_{t}$; otherwise $e_{t,i_3}\in E(M_{\tilde{k}+1}\cup\dots\cup M_{\frac{k}{2}})$ and $e_{t,i_3}\in E(G_{(t-1)n+i_3})$, without loss of generality, we could still assume that $G_{\frac{(t-1)n}{2}+i_3}\in \mathcal{G}_{t}$. If $\ell_1\neq t$ or $\ell_2\neq t$, do similar replacements as in Claim \ref{claimeven}. Hence, assume $\ell_1=\ell_2=t$ in the remaining proof. Notice that $w,w'\neq u_t$. If $v,v'\neq u_t$, replace $vv'\in E(G_{(s-1)n+i_2})$ by $vw\in E(G_{(s-1)n+i_2})$ and $ww'$ by $v'w'$ so that $e_{r,i_1}=vv'$, $e_{s,i_2}=vw$ and $e_{t,i_3}=v'w'$. Otherwise, assume $v=u_t$. Recall that for any ${n}$ graphs in the graph set $\mathcal{G}^{**}$, there are two graphs, say $G_a\in \mathcal{G}_{\ell_a}$ and $G_b\in \mathcal{G}_{\ell_b}$, such that $\ell_a\neq \ell_b$. It implies that in $M_s$, there is an edge $e_{s,i_4}=xx'\in E(G_{(s-1)n+i_4})$ such that $G_{(s-1)n+i_4}\in \mathcal{G}_{\ell_4}$ and $\ell_4\neq \ell_2$ (so $\ell_4\neq t$). Since $x\neq u_{\ell_4}$ or $x'\neq u_{\ell_4}$, we set $x'\neq u_{\ell_4}$.
    
       \noindent{\underline{{\textbf{Subcase 2.1. } $ x\neq u_t$ and $x'\neq u_t$. }}} 

       If $v'\neq  x$, replace $u_tv'\in E(G_{(s-1)n+i_2}) $ by $xv' \in E(G_{(s-1)n+i_2})$ and $xx'$ by $u_tx'$, so that $e_{r,i_1}=u_tv'$, $e_{s,i_2}=xv'$, $e_{s,i_4}=u_tx'$, and $e_{t,i_3}=ww'$. For $v'=x$, replace $u_tv'\in E(G_{(s-1)n+i_2}) $ by $wv' \in E(G_{(s-1)n+i_2})$, $v'x'$ by $u_tx'$, and $ww'$ by $v'w'$, so that $e_{r,i_1}=u_tv'$, $e_{s,i_2}=wv'$, $e_{s,i_4}=u_tx'$, and $e_{t,i_3}=v'w'$.

\noindent{\underline{{\textbf{Subcase 2.2. } $ x= u_t$ or $x'=u_t$. }}}  

 Without loss of generality, let $x'=u_t$ and $w\neq u_{\ell_4}$. Then $ x\neq u_t$. If $x\in B_{(s-1)n+i_2}$, replace $u_tv'\in E(G_{(s-1)n+i_2}) $ by $u_t x\in E(G_{(s-1)n+i_2})$, $xu_t$ by $wu_t$, and $ww'$ by $v'w'$, such that $e_{r,i_1}=u_tv'$, $e_{s,i_2}=u_{t}x$, $e_{s,i_4}=wu_t$, and $e_{t,i_3}=v'w'$. If $x\notin B_{(s-1)n+i_2}$, there exists a vertex $y\in B_{(s-1)n+i_2}$ such that $u_ty\notin E(F^*\cup M_{\tilde{k}+1}\cup\dots\cup M_{\frac{k}{2}})$. Since $d_{M_s}(y)=2$, we can find an edge $yy'\in E(M_s)$ such that $y'\neq v'$. Assume that $e_{s,i_5}=yy'\in E(G_{(s-1)n+i_5})$ and $G_{(s-1)n+i_5}\in \mathcal{G}_{\ell_5}$. Note that $y,y'\notin \{u_t,v'\}$. If $\ell_5\neq t$ and $y=u_{\ell_5}$, then $y'\neq u_{\ell_5}$ and $v'\neq u_{\ell_5}$. Thus, replace $u_tv'\in E(G_{(s-1)n+i_2})$ by $u_ty\in E(G_{(s-1)n+i_2})$ and $yy'$ by $v'y'$, so that $e_{r,i_1}=u_tv'$, $e_{s,i_2}=u_{t}y$, $e_{s,i_4}=xu_t$, $e_{s,i_5}=v'y'$, and $e_{t,i_3}=ww'$. If $\ell_5\neq t$ and $y\neq u_{\ell_5}$, then replace $u_tv'\in E(G_{(s-1)n+i_2})$ by $v'y'\in E(G_{(s-1)n+i_2})$ and $yy'$ by $u_ty$, so that $e_{r,i_1}=u_tv'$, $e_{s,i_2}=v'y'$, $e_{s,i_4}=xu_t$, $e_{s,i_5}=u_t y$, and $e_{t,i_3}=ww'$. If $\ell_5=t$, $u_tv'\in E(G_{(s-1)n+i_2}) $ by $u_{t}y \in E(G_{(s-1)n+i_2})$ and $yy'$ by $v'y'$, so that $e_{r,i_1}=u_tv'$, $e_{s,i_2}=u_{t}y$, $e_{s,i_4}=xu_t$, $e_{s,i_5}=v'y'$ and $e_{t,i_3}=ww'$. 

Notice that Subcases 2.1 and 2.2 may lead to multiple edges in the graph $F^*\cup M_{\tilde{k}+1}\cup\cdots\cup M_{\frac{k}{2}}$. However, by repeatedly applying the above steps, we will eventually obtain $E(M_i)\cap E(M_j)=\emptyset$ and $E(F^*)\cap E(M_j)=\emptyset$ for any $i,j\in \{\tilde{k}+1,\dots,\frac{k}{2}\}$. It follows that $\mathcal{G}$ admits a rainbow $k$-factor. Hence, the proof is complete. 
\end{proof}

    \begin{lemma}\label{spectralradius}
        Let $n$, $k$ and $p$ be three positive integers. If $k\geq 2$, $n\geq 4k-1$ and $1\leq p\leq \lceil\frac{n-k}{2}\rceil-1$, then $$\rho (K_{k+p-1}\vee \left((p+1)K_1\cup K_{n-k-2p}\right))< \rho(H_{n,k}).$$
    \end{lemma}

    \begin{proof}
    For  simplicity, assume that $G=K_{k+p-1}\vee \left((p+1)K_1\cup K_{n-k-2p}\right)$ and $G'=H_{n,k}=K_{k-1}\vee(K_1\cup K_{n-k})$. Partition the vertex set of $G$ as $V(G)=V((p+1)K_1)\cup V(K_{k+p-1})\cup V(K_{n-k-2p})$, where $V((p+1)K_1)=\{u_1,u_2,\ldots,u_{p+1}\}$, $V(K_{k+p-1})=\{v_1,v_2,\ldots,v_{k+p-1}\}$ and $V(K_{n-k-2p})=\{v_{k+p},v_{k+p+1},\ldots,v_{n-p-1}\}$. One can find that 
    \begin{equation*}
       G'=G-\sum_{i=k}^{k+p-1}u_1v_{i}+\sum_{i=2}^{p+1}\sum_{j=k+p}^{n-p-1}u_iv_j+\sum_{i=2}^{p}\sum_{j=i+1}^{p+1}u_iu_j.
    \end{equation*}
       
        Let $x$ be the Perron vector of $A(G)$ and let $\rho=\rho(G)$. By symmetry, the components of $x$ corresponding to the vertices in $V((p+1)K_1)$ (resp. $V(K_{k+p-1})$ or $V(K_{n-k-2p})$) are equal, denoted by $x_1$ (resp. $x_2$ or $x_3$). Then from $A(G)x=\rho x$, we have
\begin{equation}\label{eq1}
            x_2=\frac{(\rho-(n-k-2p-1))x_3}{k+p-1}.
        \end{equation}

        Let $y$ be the Perron vector of $A(G')$ and let $\rho'=\rho(G')$. Similarly, the components of $y$ corresponding to the vertices in $V(K_1)$ (resp. $V(K_{k-1})$ or $V(K_{n-k})$) are equal, denoted by $y_1$ (resp. $y_2$ or $y_3$). Then from $A(G')y=\rho'y$, we have 
        \begin{align*}
            \rho'y_1&=(k-1)y_2,\\
            \rho'y_3&=(k-1)y_2+(n-k-1)y_3.
        \end{align*}
        
        Since $G'$ contains $K_{n-1}$ as a proper subgraph, we have $\rho'\geq n-2$. Combining with the above two equations, we get 
\begin{equation}\label{eq2}
    y_3=\frac{\rho'y_1}{\rho'-(n-k-1)}.
        \end{equation}
        It is clearly that neither $G$ nor $G'$ is a complete graph. So $\rho<n-1$ and $\rho'<n-1$. Combining with \eqref{eq1} and \eqref{eq2}, we yield
\begin{align*}
    y^\top (\rho'-\rho)x  & =\rho'y^\top x-\rho y^\top x\\
             & =y^\top A(G')x-y^\top A(G)x\\
             & = \sum_{i=2}^{p+1}\sum_{j=k+p}^{n-p-1} {(x_{u_i}y_{v_j}+y_{u_i}x_{v_j})}+\sum_{i=2}^{p+1}\sum_{j=2,j\neq i}^{p+1}{x_{u_i}y_{u_j}}-\sum_{i=k}^{k+p-1}{(x_{u_1}y_{v_i}+y_{u_1}x_{v_i})}\\
             & =p(n-k-2p)(x_1y_3+x_3y_3)+p(p-1)x_1y_3-p(x_1y_3+y_1x_2)\\
             & =p[(n-k-2p+p-1-1)x_1y_3+(n-k-2p)x_3y_3-y_1x_2]\\
             & >p[(n-k-2p)x_3y_3-y_1x_2]\text{ (since } 1\leq p\leq \lceil \frac{n-k}{2}\rceil -1\text{, } n\geq 4k-1 \text{ and } k\geq 2\text{)}\\
             & = px_3y_1[\frac{\rho'(n-k-2p)}{\rho'-(n-k-1)}-\frac{\rho-(n-k-2p-1)}{k+p-1}]\text{ (by \eqref{eq1} and \eqref{eq2})}\\
             & > px_3y_1(\frac{(n-1)(n-k-2p)}{k}-\frac{k+2p}{k+p-1})\text{ (since }\rho<n-1\text{ and }\rho'<n-1\text{)}\\
             & = \frac{px_3y_1}{k(k+p-1)}((n-1)(n-k-2p)(k+p-1)-k(k+2p))\\
             & > \frac{px_3y_1}{k(k+p-1)}(2k(k+p-1)-k(k+2p)) \text{ (since }n-1>2k\text{ and } n-k-2p\geq 1\text{)}\\
             &\geq  0\text{ (since } k\geq 2\text{)}.
\end{align*}
        It implies $\rho<\rho'$. Therefore, the proof is complete.
    \end{proof}

    Now we shall give the proof of Theorem \ref{rainbowkfactor}. 
    
    \noindent{{\bf Proof of Theorem \ref{rainbowkfactor}.}} For $k=1$, the result follows from Theorem \ref{rainbowmatching}. Considering $k\geq 2$, suppose to the contrary that $\mathcal{G}=\{G_1,\ldots,G_{\frac{kn}{2}}\}$ does not admit a rainbow $k$-factor. According to Lemma \ref{structural}, the graph set $KO(\mathcal{G})=\{KO(G_1),\ldots,KO(G_{\frac{kn}{2}})\}$ does not admit a rainbow $k$-factor. For simplicity, denote $G'_t := KO(G_t)$ for $t\in [\frac{kn}{2}]$ and $\mathcal{G}':=KO(\mathcal{G})$. With Lemma \ref{boundKO} and the spectral conditions of Theorem \ref{rainbowkfactor}, we obtain that for each integer $t\in [\frac{kn}{2}]$, 
    \begin{equation}\label{eqrho}
        \rho(G'_t)\geq \rho(G_t)\geq\rho(H_{n,k}).
    \end{equation} 

    Next, we show that $\rho(G_t)=\rho(H_{n,k})$ and $ G'_t\cong H_{n,k}$ for $t\in [\frac{kn}{2}]$. Denote by $\{i,j\}$ the edge incident with the vertices $i$ and $j$.
      \setcounter{claim}{0}
      \begin{claim}\label{claim-edge-k-1}
        For each $G'_t\in\mathcal{G}'$, $\{i,j\}\in E(G'_t)$ where $i\in [k-1]$ and $j\in [n]$.
      \end{claim}
      \begin{proof}
         First of all, assert that $\{i,n\}\in E(G'_t)$ for $i\in [k-1]$ and $t \in [\frac{kn}{2}]$. Otherwise, $G'_t$ will be a proper subgraph of $H_{n,k}$. By Lemma \ref{lemsubgraph}, we have $\rho(G'_i)<\rho(H_{n,k})$, a contradiction. Then using Proposition \ref{prop}, we deduce that $\{i,j\}\in E(G'_t)$ for each $t\in [\frac{kn}{2}]$ when $i \in [k-1]$ and $j\in [n]$.  
      \end{proof}
 
    \begin{claim}\label{claim-edge-ki}
        For each $G'_t\in\mathcal{G}'$, $\{k+i,n-i\}\in E(G'_t)$ where $i\in [\lceil \frac{n-k}{2} \rceil-1 ]$.
    \end{claim}
    \begin{proof}
        Assume that there exist integers $p\in [\lceil \frac{n-k}{2} \rceil-1]$ and $t_1\in [\frac{kn}{2}]$ satisfying $\{k+p,n-p\}\notin E(G'_{t_1})$. By Proposition \ref{prop}, every edge $\{i,j\}$ of $G'_{t_1}$, where $i<j$, satisfies $i<k+p$ or $j<n-p$. Then $G'_{t_1}$ is a subgraph of $K_{k+p-1}\vee (K_{n-k-2p}\cup (p+1)K_1)$. From Lemmas \ref{lemsubgraph} and \ref{spectralradius}, we obtain that $$\rho(G'_{t_1})\leq\rho(K_{k+p-1}\vee (K_{n-k-2p}\cup (p+1)K_1))<\rho(H_{n,k}),$$ a contradiction. Thus, $\{k+i,n-i\}\in E(G'_t)$ for $t=1,2,\dots,\frac{kn}{2}$ and $i=1,2,\dots,\lceil\frac{n-k}{2}\rceil-1$.
    \end{proof}
 \begin{claim}\label{claim-edge-kn}
   For any $G'_t\in\mathcal{G}'$, $\{k,n\}\notin E(G'_t)$.
    \end{claim}
\begin{proof}
    Suppose by way of contradiction that there exists an integer $t_2\in [\frac{kn}{2}]$ such that $\{k,n\}\in E(G'_{t_2})$. Without loss of generality, assume $t_2=\frac{(k-1)n}{2}+k$ for even $n$ and $t_2=\frac{(k-1)(n-1)}{2}+k-1$ for odd $n$. Let $F$ be a $k$-regular graph with vertex set $V(F)=V$ and edge set $E(F)$. Moreover, when $n$ is even, define $E(F)=\{e_{i,j}|1\leq i\leq k,1\leq j\leq \frac{n}{2}\}$ where
    \begin{equation*}
       e_{i,j}=\left\{\begin{aligned}
       &\{j,\frac{n}{2}+i-j\}, \text{ when $1\leq j\leq  i-1,$}\\
    &\{j,n+i-j\}, \text{ when $ i\leq j\leq \frac{n}{2}$.}
        \end{aligned}
        \right.
    \end{equation*}
Claims \ref{claim-edge-k-1} and \ref{claim-edge-ki} imply that when $i\neq k$ or $j\neq k$, $e_{i,j}\in E(G'_t)$ for each $G'_t\in \mathcal{G}'$. In addition, $e_{k,k}\in E(G'_{\frac{(k-1)n}{2}+k})$ by assumption. It implies that $e_{i,j}\in E(G'_{\frac{(i-1)n}{2}+j})$ for each $i\in [k]$ and $j\in [\frac{n}{2}]$. When $n$ is odd, define $E(F)=\{e_{0,j}|1\leq j\leq \frac{k}{2}\}\cup \{e_{i,j}|1\leq i\leq k,1\leq j\leq \frac{n-1}{2}\}$ where $e_{0,1}=\{1,\frac{n+1}{2}\}$ and 
     \begin{equation*}
       e_{i,j}=\left\{\begin{aligned}
       &\{2j-2,2j-1\}, \text{ when $i=0$ and $2\leq j\leq \frac{k}{2},$}\\
      & \{j,n+1-j\}, \text{ when $i=1$ and $1\leq j\leq \frac{n-1}{2},$}\\
    & \{j+1,n+i-1-j\}, \text{ when $i\geq 2$ and $i-1\leq j\leq \frac{n-1}{2},$}\\
      & \{j,\frac{n+1}{2}+i-1-j\}, \text{ when $i\geq 3$ and $1\leq j\leq i-2.$}
        \end{aligned}
        \right.
    \end{equation*}
 Similarly, we can find that $e_{i,j}\in E(G'_{\frac{(i-1)(n-1)}{2}+j})$ for each $i\in [k]$ and $j\in [\frac{n-1}{2}]$ and $e_{0,j}\in E(G'_{\frac{k(n-1)}{2}+j})$ for $j\in [\frac{k}{2}]$. Thus, the graph $F$ is a rainbow $k$-factor of $\mathcal{G}'$, a contradiction. 
\end{proof}

    According to Claim \ref{claim-edge-kn} and Proposition \ref{prop}, each $G'_t\in \mathcal{G}'$ is a subgraph of $H_{n,k}$. By Lemma \ref{lemsubgraph}, we deduce that for $t=1,2,\dots,\frac{kn}{2}$, $\rho(G'_t)\leq\rho(H_{n,k})$, with equality if and only if $G'_t\cong H_{n,k}$. Combining with \eqref{eqrho}, we obtain $\rho(G'_t)=\rho(G_t)=\rho(H_{n,k})$ and $G'_t\cong H_{n,k}$ for each $t=1,2,\ldots,\frac{kn}{2}$. By Corollary \ref{cor} and Lemma \ref{factlemma}, it follows that $G_1=\dots=G_{\frac{kn}{2}} \cong H_{n,k}$. Therefore, the proof is complete.\qed


\begin{thebibliography}{99}
    \bibitem{ADGM20}
    R. Aharoni, M. DeVos, S. Gonz\'{a}lez Hermosillo de la Maza, A. Montejano, and R. \v{S}\'{a}mal, A rainbow version of Mantel’s theorem, \textit{Adv. Comb.}, Paper No. 2, 12pp, 2020.

    \bibitem{BM76}
    J.~A. Bondy and U.~S.~R. Murty, Graph theory with applications, American Elsevier Publishing Co., Inc., New York, 1976.
    

    \bibitem{B21}
    P. Bradshaw, Transversals and bipancyclicity in bipartite graph families, \textit{Electron. J. Combin.}, 28(4): Paper No. 4.25, 20pp, 2021.


    \bibitem{BH2005}
    A.~E. Brouwer and W.~H. Haemers, Eigenvalues and perfect matchings, \textit{Linear Algebra Appl.}, 395: 155--162, 2005. 
    
    \bibitem{BH12}
    A.~E. Brouwer and W.~H. Haemers, Spectra of graphs, Universitext, Springer, New York, 2012.


    \bibitem{CSWW24+}
    Y. Cheng, W. Sun, G. Wang, and L. Wei, Transversal Hamilton paths and cycles, arXiv:2406.13998.

    \bibitem{CHOP2021}
    E.-K. Cho, J.~Y. Hyun, S. O and J.~R. Park, Sharp conditions for the existence of an even $[a,b]$-factor in a graph, \textit{Bull. Korean Math. Soc.}, 58(1): 31--46, 2021.
    
    \bibitem{Cioaba2005}
    S.~M. Cioab\u a, Perfect matchings, eigenvalues and expansion, \textit{C. R. Math. Acad. Sci. Soc. R. Can.}, 27(4): 101--104, 2005. 
  
    \bibitem{CG2007}
    S.~M. Cioab\u a{} and D.~A. Gregory, Large matchings from eigenvalues, \textit{Linear Algebra Appl.}, 422(1): 308--317, 2007.

    
    \bibitem{CGH2009}
    S.~M. Cioab\u a, D.~A. Gregory and W.~H. Haemers, Matchings in regular graphs from eigenvalues, \textit{J. Combin. Theory Ser. B}, 99(2): 287--297,  2009.

    \bibitem{C09}
    P. Csikv\'ari, On a conjecture of V. Nikiforov, \textit{Discrete Math}., 309(13): 4522--4526, 2009.
    
    \bibitem{FLL22}
    D. Fan, H. Lin and H. Lu, Spectral radius and $[a,b]$-factors in graphs, \textit{Discrete Math.}, 345(7): Paper No. 112892, 9pp, 2022.

    \bibitem{FLLO23}
    D. Fan, H. Lin, H. Lu and S. O, Eigenvalues and factors: a survey, arXiv:2312.15902.
    
    \bibitem{FN10}
    M. Fiedler and V. Nikiforov, Spectral radius and Hamiltonicity of graphs, \textit{Linear Algebra Appl.}, 432(9): 2170--2173, 2010.
    
    \bibitem{GLMM23}
    M. Guo, H. Lu, X. Ma, and X. Ma, Spectral radius and rainbow matchings of graphs, \textit{Linear Algebra Appl.}, 679: 30--37, 2023.

    \bibitem{GHMPS23}
    P. Gupta, F. Hamann, A. M\"{u}yesser, O. Parczyk, and A. Sgueglia, A general approach to transversal versions of Dirac-type theorems, \textit{Bull. Lond. Math. Soc.}, 55(6): 2817–2839, 2023.

    \bibitem{HFGLSTKZ24}
    Z. He, P. Frankl, E. Gy\H{o}ri, Z. Lv, N. Salia, C. Tompkins, K. Varga, and X. Zhu, Extremal results for graphs avoiding a rainbow subgraph, \textit{Electron. J. Combin.}, 31(1): Paper No. 1.28, 11pp, 2024.

    \bibitem{HLF24}
    X. He, Y. Li and L.~H. Feng, Spectral radius and rainbow Hamilton paths of a graph, \textit{Discrete Math.}, 347(10): Paper No. 114128, 11pp, 2024.

    
    \bibitem{JK20}
    F. Joos and J. Kim, On a rainbow version of Dirac's theorem, \textit{Bull. Lond. Math. Soc.}, 52(3): 498--504, 2020.

    \bibitem{KSSV04}
    P. Keevash, M. Saks, B. Sudakov, and J. Verstra\"{e}te, Multicolour Tur\'{a}n problems, \textit{Adv. in Appl. Math.}, 33(2): 238--262, 2004.

\bibitem{Kelmans1981}
 A.~K. Kelmans, On graphs with randomly deleted edges, \textit{Acta Math. Acad. Sci. Hungar.}, 37(1-3): 77–88, 1981.
 
    \bibitem{LLL23}
    L. Li, P. Li, and X. Li, Rainbow structures in a collection of graphs with degree conditions, \textit{J. Graph Theory}, 104(2): 341–359, 2023.


\bibitem{LLD19}
    M. Liu, H-J, Lai, K.~Ch. Das, Spectral results on Hamiltonian problem, \textit{Discrete Math.}, 342(6): 1718–1730, 2019.

    \bibitem{Lu2010}
    H. Lu, Regular factors of regular graphs from eigenvalues, \textit{Electron. J. Combin.}, 17(1): Research Paper 159, 12pp, 2010.
    
    \bibitem{Lu2012}
    H. Lu, Regular graphs, eigenvalues and regular factors, \textit{J. Graph Theory}, 69(4): 349--355, 2012.
 
    \bibitem{M15}
    C. Magnant, Density of Gallai multigraphs,\textit{ Electron. J. Combin.}, 22(1): Paper 1.28, 6pp, 2015.

     \bibitem{O2021}
     S. O, Spectral radius and matchings in graphs, \textit{Linear Algebra Appl.}, 614: 316--324, 2021.

    \bibitem{SWW24+}
    W. Sun, G. Wang and L. Wei, Transversal Structures in Graph Systems: A Survey, arXiv:2412.01121.
    
    \bibitem{ZV24+}
    Y. Zhang and E. R. van Dam, Rainbow Hamiltonicity and the spectral radius, \textit{Discrete Math.}, 348(11): Paper No. 114600, 7 pp, 2025.

    \end{thebibliography}
\end{document}